\documentclass[12pt]{article}
\usepackage[all]{xy}
\title{s-pure extensions of locally compact abelian groups}
\author{H.Sahleh\\Department of Mathematics, University of Guilan, P.O.BOX 1914\\Rasht-Iran\\
e-mail: sahleh@guilan.ac.ir\\
A.A. alijani\\ Department of Mathematics, University of Guilan\\e-mail: taleshalijan@phd.guilan.ac.ir }
\begin{document}
\date{}
\maketitle
\begin{abstract}
\newcommand{\stk}{\stackrel}
 A subgroup $H$ of a locally compact abelian (LCA) group $G$ is called s-pure if $ \overline{H \bigcap nG}=H $ for every positive integer $n$. A proper short exact sequence $0\to A\stackrel{\phi}{\to} B\to C\to 0$ in the category of LCA groups is said to be {\it s-pure} if $\phi(A)$ is an s-pure subgroup of $G$. We establish conditions under which the s-pure exact sequences split and determine those LCA groups which are s-pure injective. We also gives a necessary condition for an LCA group to be s-pure projective in $\pounds$.

AMS.subj.class: 22B05
\end{abstract}

\centerline {Introduction}
\newcommand{\stk}[1]{\stackrel{#1}{\longrightarrow}}
  All groups considered in this paper are Hausdorff topological abelian groups and they will be written additively. For a group $G$ and a positive integer $n$, we denote by $nG$, the subgroup of $G$ defined by $nG=\{nx:x\in G\}$ and $G[n]$, the subgroup of $G$ defined by $G[n]=\{x\in G;nx=0\}$. In a multiplicative group, we will use $G^{n}$ instead of $nG$ and define $G^{n}=\{x^{n}:x\in G\}$. Let $\pounds$ denote the category of locally compact abelian ($LCA$) groups with continuous homomorphisms as morphisms. A morphism is called proper if it is open onto its image, and a short exact sequence $0 \to A \stk{\phi} B \stk{\psi} C \to 0$ in $\pounds$ is said to be an extension of $A$ by $C$ if $\phi$ and $\psi$ are proper morphism. We let $Ext(C,A)$ denote the group of extensions of $A$ by $C$ [6]. Let $\overline{S}$ denotes the closure of $S\subseteq G$. We say that a closed subgroup $H$ of an LCA group $G$ is s-pure if $ \overline{H \bigcap nG}=H $ for every positive integer $n$. A subgroup $H$ of a group $G$ is said to be pure if $H \bigcap nG=nH $ for every positive integer $n$ [3]. A pure subgroup need not be s-pure and vice versa (Example 1.9). In Section 1, we show that an s-pure subgroup is pure if and only if it is densely divisible (Lemma 1.10). An LCA group $G$ is said to be pure simple if $G$ contains no nontrivial closed pure subgroup [1]. Armacost [1] has determined the pure simple LCA group $G$. Also, Armacost has determined the LCA group $G$ such that every closed subgroup of $G$ is pure [1]. We say that an LCA group $G$ is s-pure simple if $G$ contains no nonzero s-pure subgroup. We say that a LCA group $G$ is s-pure full if every closed subgroup of $G$ is s-pure. We show that a LCA group $G$ is s-pure full if and only if it is divisible (Theorem 1.15). Also, we show that a compact group $G$ is s-pure simple if and only if it is totally disconnected(Theorem 1.16). A proper short exact sequence  $0 \to A \stk{\phi} B\to C \to0$ in $\pounds$ is said to be s-pure if $\phi(A)$ is s-pure in $B$. In section 2, we study s-pure exact sequence in $\pounds$. In [4], Fulp studied pure injective and pure projective in $\pounds$. In section 3, we study s-pure injective and s-pure projective in $\pounds$. An LCA group $G$ is an s-pure injective group in $\pounds$ if and only if $ G \cong R^{n}\bigoplus (R/Z)^{\sigma}$ (Theorem 3.2). If $G$ is an s-pure projective group in $\pounds$ then $G\cong R^{n}\bigoplus G' $ where $G'$ is a discrete torsion-free, non divisible group (Theorem 3.4).\\

   The additive topological group of real numbers is denoted by $R$, $Q$ is the group of rationals with the discrete topology and $Z$ is the group of integers. Also, $Z(n)$ is the cyclic group of order $n$ and $Z(p^{\infty})$ denotes the quasicyclic group. For any group $G$, $G_0$ is the identity component of $G$, $tG$ is the maximal torsion subgroup of $G$ and $1_{G}$ is the identity map $G\to G$. An element $g\in G$ is called compact if the smallest closed subgroup which its contains is compact [8, Definition 9.9]. We denote by $bG$, the subgroup of all compact elements of $G$. If $\{G_{i}\}_{i\in I}$ is a family of groups in $\pounds$, then we denote their direct product by $\prod_{i\in I}G_{i}$. If all the $G_{i}$ are equal, we will write $G^{I}$ instead of $\prod_{i\in I}G_{i}$. For any group $G$ and $H$, $Hom(G,H)$ is the group of all continuous homomorphisms from $G$ to $H$, endowed with the compact-open topology. The dual group of $G$ is $\hat{G}=Hom(G,R/Z)$ and $(\hat{G},S)$ denotes the annihilator of $S\subseteq G$ in $\hat{G}$. For a group $G$, we define $G^{(1)}=\bigcap_{n=1}^{\infty} nG$.

\section{s-pure subgroups}

Let $G\in \pounds$. In this section, we introduce the concept and study some properties of an s-pure subgroup of $G$.\\

{\bf 1.1. Definition}. {\it A closed subgroup $H$ of a group $G$ is called s-pure if $ \overline{H \bigcap nG}=H $ for every positive integer $n$.}\\

{\bf 1.2. Note}.
\begin{itemize}
\item [(a)] A closed divisible subgroup of a group is s-pure.
\item [(b)] A closed subgroup of a divisible group is s-pure.\\
\end{itemize}

{\bf 1.3. Remark}. {\it Let $G\in\pounds$. Then $G$ has two trivial subgroups,$\{0\}$ and $G$. Clearly, $\{0\}$ is s-pure. But $G$ need not be an s-pure in itself.\\

Recall that a group $G$ is said to be densely divisible if it has a dense divisible subgroup.\\

{\bf 1.4. Lemma}. {\it A group $G$ is densely divisible if and only if $\overline{nG}=G$ for every positive integer $n$.}\\

Proof. See [2, 4.16(a)].\\

{\bf 1.5. Corollary}. {\it Let $G\in \pounds$. Then, $G$ is s-pure in itself if and only if $G$ is densely divisible.}\\

Proof. It is clear by Lemma 1.4.\\

{\bf 1.6. Lemma}. {\it Let $G\in \pounds$. Then, $\overline{G^{(1)}}$ is an s-pure subgroup of $G$.}\\

Proof. It is clear that $G^{(1)}\subseteq \overline{G^{(1)}} \bigcap mG$ for every positive integer $m$. So, $\overline{G^{(1)}}\subseteq\overline{\overline{G^{(1)}} \bigcap mG}\subseteq\overline{G^{(1)}}$ for all $m$. Hence, $\overline{G^{(1)}}$ is an s-pure subgroup.\\

{\bf 1.7. Remark}. Let $G\in \pounds$ and $H$ be an s-pure subgroup of $G$. Then, $H\subseteq \overline{nG}$ for all positive integers $n$. Hence, $H\subseteq \bigcap_{n=1}^{\infty} \overline{nG}=(G,t\hat{G})$ [8].\\

 Now, we present an example of a LCA group $G$ and a closed subgroup $H$ of $G$ such that $H\subseteq (G,t\hat{G})$, but $H$ is not an s-pure subgroup.\\

{\bf 1.8. Example}. Let $S^{1}$ be the (multiplicative) circle group of unitary complex numbers and $\sigma$ any infinite cardinal number. Let $G$ be the subgroup of $(S^{1})^{\sigma}$ consisting of all $(x_{\iota})$ such that $x_{\iota}=\pm 1$ for all but a finite number of $\iota$. Let $K$ be the subgroup of $G$ consisting of all $(x_{\iota})$ such that $x_{\iota}=1$ for all but a finite number of $\iota$. By [8, section 24.44(a)], $G$ is a locally compact abelian group, and $\hat{G}$ is torsion-free. Let $H=\{(x)_{\iota},(y)_{\iota}\}$ where $x_{\iota}=1$ and $y_{\iota}=-1$ for $\iota\neq\iota_{1},...,\iota_{m}$ and $x_{\iota}=y_{\iota}=0$ for $\iota= \iota_{1},...,\iota_{m}$. Then, $H$ is a closed subgroup of $G$, and $H\subseteq G=(G,t\hat{G})$. Now, suppose that $n$ is even. Then, $\overline{H\bigcap G^{n}}=\overline{H\bigcap K}=\{(x)_{\iota}\}$. Hence, $H$ is not s-pure.\\

Recall that a subgroup $H$ of a group $G$ is called pure if $nH=H\bigcap nG$ for every positive integer $n$[3]. A pure subgroup need not be s-pure, and an s-pure subgroup need not be pure.\\

{\bf 1.9. Example}. Since $R$ is divisible, so the subgroup $Z$ of $R$ is s-pure. But it is not a pure subgroup. Let $p$ be a prime and $G=\prod_{n=1}^{\infty}Z(p^{n})$, with discrete topology. Then, $tG$ is a pure subgroup of $G$. Since $(1,0,0,...)\in tG$ and $(1,0,0,...)\notin p(tG)$, so it is not s-pure.\\

 {\bf 1.10. Lemma}. {\it A pure subgroup is s-pure if and only if it is densely divisible.}\\

 Proof. Let $H$ be a pure subgroup of $G$. If $H$ is an s-pure subgroup, then $\overline{nH}=H$ for every positive integer $n$. So, by Lemma 1.4, $H$ is densely divisible. Conversely, let $H$ be a densely divisible, pure subgroup of $G$. Then, $ \overline{H \bigcap nG}=\overline{nH} $ for every positive integer $n$. By Lemma 1.4, $\overline{nH}=H$ for all $n$. So, $ \overline{H \bigcap nG}=H $ for all $n$. Hence, $H$ is an s-pure subgroup in $G$. \\

Let $G$ be a group in $\pounds$. Then $G$ is called s-pure simple if $G$ contains no nonzero s-pure subgroups. Similarly, $G$ is called s-pure full if every closed subgroup of $G$ is s-pure.\\

{\bf 1.11. Lemma}. {\it Let $G_{1}$ and $G_{2}$ be two groups in $\pounds$. If $G_{1}\times G_{2}$ is s-pure full, then $G_{1}$ and $G_{2}$ are s-pure full.}\\

Proof. Let $G_{1}$,$G_{2}\in \pounds$ and $H$ be a closed subgroup of $G_{1}$. Then, $H\times G_{2}$ is a closed subgroup of $G_{1}\times G_{2}$. So, $\overline{(H\times G_{2})\bigcap (nG_{1}\times nG_{2})}=H\times G_{2}$ for any positive integer $n$. Hence, $\overline{(H\bigcap nG_{1})\times (G_{2}\bigcap nG_{2})}=H\times G_{2}$. Therefore, $\pi_{1}(\overline{(H\bigcap nG_{1})\times (nG_{2})})=\pi_{1}(H\times G_{2})$ where $\pi_{1}$ is the first projection map of $G_{1}\times G_{2}$ onto $G_{1}$. Consequently, $\overline{H\bigcap nG_{1}}=H$. Similarly, it can be show that $G_{2}$ is s-pure full. \\

{\bf 1.12. Remark}. {\it Recall that a discrete group is densely divisible if and only if it is divisible.} \\

{\bf 1.13. Remark}. {\it Let $G$ be a densely divisible group and $H$ a closed subgroup of $G$. Since $(\hat{G},H)$ is a subgroup of $\hat{G}$ and $\hat{G}$ is torsion-free, so $G/H$ is densely divisible.}\\

{\bf 1.14. Remark}. {\it Let $G$ be a densely divisible group and $H$ an open, pure subgroup of $G$. An easy calculation shows that $H$ is divisible.}\\

{\bf 1.15. Theorem}. {\it Let $G\in\pounds$. Then, $G$ is s-pure full if and only if $G$ is divisible.}\\

Proof. Let $G$ be an s-pure full group in $\pounds$. By [8, Theorem 24.30], $G\cong R^{n}\bigoplus G'$ , where $G'$  is an LCA group which contains a compact open subgroup. By Lemma 1.11, $G'$ is s-pure full. So, by Corollary 1.5, $G'$ is densely divisible. By Remark 1.13 , $G'/bG'$ is densely divisible. On the other hand, $G'/bG'$ is discrete and torsion-free (see the proof of Theorem 2.7 [9]). Hence, by Remark 1.12, $G'/bG'$ is divisible. By Remark 1.14, $bG'$ is divisible. Consequently, the short  exact sequence $0\to bG'\to G'\to G'/bG'\to 0$  splits. Hence, $G'\cong bG'\bigoplus G'/bG'$ and $G'$ is divisible. Therefore, $G$ is divisible. The converse is clear by Note 1.2.b. \\

{\bf 1.16. Theorem}. A compact group $G$ is an s-pure simple group if and only if it is totally disconnected.\\

Proof. Let $G$ be a compact group. If $G$ is an s-pure simple group, then by Note 1.2(a), $G_{0}=0$ because $G_{0}$ is a closed divisible subgroup of $G$. So $G$ is totally disconnected. Conversely, Let $G$ be a compact, totally disconnected group and $H$ an s-pure subgroup of $G$. By Remark 1.7, $H\subseteq(G,t\hat{G})$. Since $\hat{G}$ is a discrete and a torsion group, so $t\hat{G}=\hat{G}$. Hence, $H=0$.

\section{s-pure exact sequence}

In this section, we introduce the concept and study some properties of s-pure extensions in $\pounds$.\\

 {\bf 2.1. Definition}. {\it An extension $0 \longrightarrow A \stk{\phi} B \longrightarrow C \longrightarrow 0$ in $\pounds $ is called s-pure if $\phi(A) $ is s-pure in $B$. }\\

{\bf 2.2. Remark}. {\it Let $A$ be a divisible group in $\pounds$ and $E:0 \longrightarrow A \stk{\phi} B \longrightarrow C \longrightarrow 0$ an extension in $\pounds$. Then $\phi(A)$ is a closed divisible subgroup of $B$. So, by Note 1.2($a$),$E$ is an s-pure extension.\\

 {\bf 2.3. Lemma}. {\it Let $A,C$ be groups in $\pounds$. Then the extension $0\to A\to A\bigoplus C\to C\to 0$ is an s-pure extension if and only if $A$ is densely divisible.}\\

 Proof. The extension $0\to A\to A\bigoplus C\to C\to 0$ is pure. Hence, by Lemma 1.10, it is s-pure if and only if $A$ is densely divisible.\\

  {\bf 2.4. Remark}. Lemma 2.3 shows that the set of all s-pure extensions of $A$ by $C$ need not be a subgroup of $Ext(C,A)$.\\

 The dual of an extension $E:0 \to A \to B \to C \to 0$ is defined by $\hat{E}:0 \to \hat{C} \to \hat{B} \to \hat{A} \to 0$. The following example shows that the dual of an s-pure extension need not be s-pure.\\

{\bf 2.5. Example} There exists a non splitting extension $$ E:0\to Z(p^{\infty})\to B\to C\to 0$$ of $Z(p^{\infty})$  with compact group $C$ which is not torsion-free [2, Example 6.4]. By Note 1.2(a), $E$ is s-pure. Since $\widehat{Z(p^{\infty})}$ is torsion-free, so $\hat{E}$ is pure. By Lemma 1.10, $\hat{E}$ is s-pure if and only if $\hat{C}$ is densely divisible. But $C$ is compact. So, $\hat{C}$ is discrete. Hence, $\hat{E}$ is s-pure if and only if $\hat{C}$ is a discrete divisible group. Consequently, $\hat{E}$ is s-pure if and only if $C$ is a compact torsion-free group. Since $C$ is not torsion-free, it follows that $\hat{E}$ is not s-pure.\\

Recall that two extensions $0 \to A \stk{\phi_{1}} B \stk{\psi_{1}} C \to 0$ and $0 \to A \stk{\phi_{2}} X \stk{\psi_{2}} C \to 0$ is said to be equivalent if there is a topological isomorphism $\beta:B\to X$ such that the following diagram
\[
\xymatrix{
0 \ar[r] & A \ar^{\phi_{1}}[r] \ar^{1_{A}}[d] & B \ar^{\psi_{1}}[r] \ar^{\beta}[d] & C \ar[r] \ar^{1_{C}}[d]
& 0 \\
0 \ar[r] & A \ar^{\phi_{2}}[r] &   X \ar^{\psi_{2}}[r] & C \ar[r] & 0
}
\]
is commutative.\\

{\bf 2.6. Lemma} {\it An extension equivalent to an s-pure extension is s-pure. }\\

Proof. Suppose that $$ E_1: 0 \to A \stk{\phi_1} B \to C \to 0  ,  E_2: 0 \to A \stk{\phi_2} X \to C \to 0 $$  be two equivalent extension such that
$E_1$ is s-pure. Then there is a topological isomorphism $ \beta :B \to X $ such that $ \beta\phi_1 =\phi_2 $ . Since $E_1$ is s-pure, $ \phi_1(A) = \overline{\phi_1(A)\bigcap nB}$. Then $ \beta\phi_1(A)=\beta (\overline{\phi_1(A)\bigcap nB})$. So, $ \phi_2(A)=\overline{\phi_2(A)\bigcap nX}$. Hence, $E_2$ is s-pure.\\

{\bf 2.7. Corollary}. {\it If the s-pure extension  $ 0 \to A \to B \to C \to 0 $ splits, Then $A$ is densely divisible.}\\

Proof. Let $ 0 \to A \to B \to C \to 0 $ be a split, s-pure extension. Then, it is equivalent to $0\to A\to A\bigoplus C\to C\to 0$. So, $0\to A\to A\bigoplus C\to C\to 0$ is s-pure. Hence, by Lemma 2.3, $A$ is densely divisible.\\

{\bf 2.8. Remark}. The converse of Corollary 2.7 may not hold. Consider Example 2.5.\\

We will now show that a pullback or pushout of an s-pure extension need not be s-pure. For more on a pullback and a pushout of an extension in $\pounds$, see [6]. \\

\newcommand{\cT}{\mathcal{T}}
\newcommand{\dwn}[1]{{\scriptstyle #1}\downarrow}

{\bf 2.9. Example} Let $\alpha$ be the map $ \alpha: Z \to Z:n\longmapsto 2n $. \emph{}Consider the s-pure extension $E : 0 \to Z_{2}\to R/Z \to R/Z \to 0 $ which is the dual of $0 \to Z \stk{\alpha} Z \to Z_{2}\to 0 $. Let $f:Q \to R/Z$ be any continuous homomorphism. Since $Q$ is torsion-free, so the standard pullback of $E$ is pure, but not s-pure by Lemma 1.10 because $Z_{2}$ is not densely divisible. Now consider the s-pure extension $E':0 \to Z \to Q \to Q/Z \to 0$ . Then the map $\alpha$ induces a pushout diagram
\[
\xymatrix{
E':0 \ar[r] & Z \ar^{}[r] \ar^{\alpha}[d] & Q \ar^{}[r] \ar^{}[d] & Q/Z \ar[r] \ar^{1_{Q/Z}}[d]
& 0 \\
\alpha E':0 \ar[r] & Z \ar^{\mu}[r] &   (Z\bigoplus Q)/ H \ar^{}[r] & Q/Z \ar[r] & 0
}
\]

Where $H= \{(2n,-n);n\in Z\}$ and $ \mu :n \longmapsto (n,0)+H$. If $ \alpha E'$ is s-pure, then $ \mu (Z)\subseteq 2((Z \bigoplus Q)/H)$ which is a contradiction.\\

\section{s-pure injectives and s-pure projectives}

In this section, we define the concept of s-pure injective and s-pure projective in $\pounds$ and express some of their properties .\\

{\bf 3.1. Definition} {\it Let $G$ be a group in $\pounds$. We call $G$ an s-pure injective group in $\pounds$ if for every s-pure exact sequence $$ 0 \to A \stk{\phi} B \to C \to 0 $$ and continuous homomorphism $ f:A\to G$, there is a continuous homomorphism $ \bar{f}:B\longrightarrow G$ such that $ \bar{f}\phi=f$. Similarly, we call $G$ an s-pure projective group in $\pounds$ if for every s-pure exact sequence $$ 0 \to A \to B \stk{\psi} C \to 0 $$ and continuous homomorphism $ f:G \to C$, there is a continuous homomorphism $\bar{f}:G\to B$ such that $\psi\bar{f}=f$ .}\\

{\bf 3.2. Theorem} {\it  Let $G\in\pounds$. The following statements are equivalent:
\begin{enumerate}
\item $G$ is an s-pure injective in $\pounds$.
\item $ G \cong R^{n}\bigoplus (\frac{R}{Z})^{\sigma}$ where $\sigma$ is a cardinal number.\\
\end{enumerate}

 Proof. $1\Longrightarrow 2$: Let $G$ be an s-pure injective in $\pounds$. For a group $X$ in $\pounds$, consider the s-pure extension $$ E: 0\to G \stk{\phi} B \stk{} X\to 0$$ Then there is a continuous homomorphism $\bar{\phi}:B \to G$ such that $\bar{\phi}\phi=1_G$. Consequently, $E$ splits. In particular, the s-pure extension $ 0\to G \stk{} G^{*} \to G^{*}/G\to 0$ splits where $G^{*}$ is the minimal divisible extension of $G$. Hence, $G$ is divisible. So, by Remark 2.2, every extension of $G$ by $X$ is an s-pure extension. On the other hand, every s-pure extension of $G$ by $X$ splits. Hence, $Ext(X,G)=0$. By [10, Theorem 3.2], $G \cong R^{n}\bigoplus (R/Z)^{\sigma}$.

 $2\Longrightarrow 1$: It is clear.\\

Recall that a discrete group $G$ is called reduced if it has no nontrivial divisible subgroup.\\

 {\bf 3.3. Lemma}   {\it $Q$ is not an s-pure projective group.}\\

 Proof. Consider the s-pure exact sequence $0 \to Z \stk{} R \stk{\pi} R/Z \to 0$ where $\pi$ is the natural mapping. Assume that $Q$ is an s-pure projective group and $f\in Hom(Q,R/Z)$. Then, there is $\bar{f}\in Hom(Q,R)$ such that $\pi \bar{f}=f$. Hence, $\pi^{\ast}:Hom(Q,R)\to Hom(Q,R/Z)$ is surjective. Now consider the following exact sequence $$0\to Hom(Q,Z)\to Hom(Q,R)\stk{\pi_{\ast}} Hom(Q,R/Z)\to Ext(Q,Z)\to Ext(Q,R)$$ Since $Q$ is divisible and $Z$ is reduced, so $Hom(Q,Z)=0$. Hence,  $\pi^{\ast}$ is one to one. This shows that  $\pi^{\ast}$ is an isomorphism. On the other hand, $Ext(Q,R)=0$. Consequently, $Ext(Q,Z)=0$ which is a contradiction.\\

{\bf 3.4. Theorem} {\it Let $G\in\pounds$. If $G$ is an s-pure projective in $\pounds$, then $G\cong R^{n}\bigoplus G' $ where $G'$ is a discrete torsion-free , reduced group.}\\

 Proof. It is known that an $LCA$ group $G$ can be written as $G \cong R^{n} \bigoplus G'$ where $G'$ contains a compact open subgroup [8, Theorem 24.30]. An easy calculation shows that if $G$ is an s-pure projective group, then $G'$ is an s-pure projective in $\pounds$. Let $f\in Hom$($G'$,$\frac{R}{Z}$). Then there exists a continuous homomorphism $\tilde{f}: G'\to R$ such that the following diagram is commutative:
\[
\xymatrix{
&   & & G' \ar_{\tilde{f}}[dl] \ar^{f}[d]
&  \\
0 \ar[r] & Z \ar^{}[r] & R \ar^{\pi}[r] & R/Z \ar[r] & 0
}
\]
Consider the following exact sequence $$0\to Hom(G',Z)\to Hom(G',R)\stk{\pi_{\ast}} Hom(G',R/Z)\to Ext(G',Z)\to 0$$ Since $\pi_{*}$ is surjective, so $Ext(G',Z)=0$. Let $K$ be a compact open subgroup of $G'$. Then the inclusion map $i:K\to G'$ induces the surjective homomorphism $i_{*}:Ext(G',Z)\to Ext(K,Z)$. So, $Ext(K,Z)=0$. Hence, $Ext(R/Z,\hat{K})=0$. By [7, Proposition 2.17], $\hat{K}=0$. So, $K=0$. Hence, $G'$ is discrete. If $G'$ contains a subgroup of the form $Z(n)$, then $Z(n)$ is a nontrivial compact open subgroup of $G'$ which is a contradiction. So $G'$ is torsion-free. Suppose $G'$ has a nontrivial divisible subgroup. Then $G'$ has a direct summand $H\cong Q$. But then $H$ is s-pure projective, contradicting Lemma 3.3. Therefore, $G'$ is reduced.\\

\end{document}